\documentclass[a4paper,11pt]{article}
\usepackage{amsfonts}
\usepackage{amsmath}
\usepackage{amssymb}
\usepackage{graphicx}
\title{Hochschild homology of certain Soergel bimodules}
\author{Marko Sto\v si\'c\\
Instituto de Sistemas e Rob\'otica and CAMGSD,\\
Instituto Superior T\'ecnico, 
Av. Rovisco Pais 1\\
1049-001 Lisbon, 
Portugal\\
e-mail: mstosic@math.ist.utl.pt
}
\date{}

\newtheorem{theorem}{Theorem}

\newtheorem{corollary}[theorem]{Corollary}

\newtheorem{remark}{Remark}

\def\max{\mathop{\rm max}}

\def\kraj{\hfill\rule{6pt}{6pt}}

\def\deg{\mathop{\rm deg}}
\def\qdim{\mathop{\rm qdim}}
\def\det{\mathop{\rm det}}

\def\C{\mathbb{C}}

\def\N{\mathbb{N}}

\begin{document}
\maketitle
\begin{abstract}
In this paper we compute Hochschild homology of certain Soergel bimodules. Moreover, we describe explicitly the graded 
bimodule maps between Soergel bimodules. This computations are motivated by the categorifications of the colored HOMFLY-PT 
polynomial for links via Hochschild homology of Soergel bimodules.
\end{abstract}

\section{Introduction}

The Soergel bimodules were introduced by Soergel in \cite{soergel1,soergel2} in the 
context of the infinite-dimensional representation theory of simple Lie algebra and 
Kazhdan-Lusztig theory. They have nice explicit expression as the tensor products of 
the rings of polynomials invariant under the action of a symmetric 
group, tensored over rings of the same form.
Moreover, there are various quite different interpretations of the Soergel bimodules 
(e.g. the geometric interpretation \cite{webster}), that make them very important.

Recently, Soergel bimodules were used by Rouquier \cite{rouq} to categorify braid 
groups. Khovanov \cite{kovhoh} extended 
this construction by adding Hochschild homology of the Soergel bimodules involved, 
to obtain a link invariant --- the 
triply-graded categorification of the HOMFLY-PT polynomial. The Hochschild homology 
works perfectly 
to give invariance under the Markov moves (Reidemeister 1 move).  Also, Khovanov has 
shown that this construction is 
isomorphic to the previous one by Khovanov and Rozansky \cite{KR2}, that uses the Koszul complexes.  

Moreover, in \cite{nas}, we have extended the Khovanov's construction for the categorification of 1,2-colored HOMFLY-PT 
polynomial.
A larger set of Soergel bimodules was used in the definition, and the computations proving the invariance under the 
Reidemeister moves 
are much more involved.  Also, we have sketched the construction in the general case - i.e. for the categorification of 
arbitrary colored HOMFLY-PT polynomial. However, because of the difficulty of the explicit expressions, the proof of the 
invariance was postponed for the subsequent paper.\\

In this paper we prove some of the results stated in \cite{nas} for 2-colored HOMFLY-PT polynomial and conjectured to hold in 
the general case. Moreover, we obtain some further properties of  Hochschild homology of Soergel bimodules. More precisely, 
we prove that the map between 
$R_{k,l}$ and $R_{k,l}\otimes_{k+l}R_{k,l}$,  given in \cite{nas} for arbitrary $k,l\in\N$, is indeed a bimodule map. 
Moreover, we prove that such a 
bimodule map is the unique one of the degree $2kl$, and that there are no nonzero graded bimodule maps of the degree strictly 
less than $2kl$. The obtained bimodule map is closely related with Khovanov-Lauda calculus of the categorifications of 
quantum groups \cite{aron}. Moreover, the particular form of this (unique) bimodule map, gives the orthogonality of a 
basis in a certain Frobenius algebra appearing in \cite{aron}.

Furthermore, we compute the Hochschild homology of the bimodule $R_{k,l}\otimes_{k+l}R_{k,l}$. In particular, 
this gives the categorification of one digon move axiom of the calculus for (arbitrary) colored HOMFLY-PT polynomial.

Although this work was motivated by the categorification of link invariants, the contents of this paper is purely algebraic. 
We use the ``Koszul picture", namely to each bimodule we correspond a certain quotient polynomial ring.
The computation of the Hochschild homology of a bimodule is reduced to the computation of the Koszul complex of the 
corresponding polynomial ring.

The polynomial rings, being invariant under the action of  symmetric groups, are generated by Schur polynomials. Their basic 
properties as well as the specific notational conventions that are used throughout the paper, are given in section \ref
{Sur-sim}. In section \ref{3}, we explain more precisely the way how we shall compute 
the Hochschild homology of bimodules. Section \ref{4} gives the explicit expression  of a bimodule map between $R_{k,l}$ and 
$R_{k,l}\otimes_{k+l}R_{k,l}$. Finally, section \ref{5} contains the computation of 
the Hochschild homology of the bimodule $R_{k,l}\otimes_{k+l}R_{k,l}$.

\section{Schur and symmetric polynomials}\label{Sur-sim}

Among various bases of homogeneous symmetric polynomials, we will concentrate on elementary symmetric polynomials and Schur 
polynomials. In this section we recall their definition and give the properties that we shall use in the paper. For more 
details see e.g. \cite{fulton}, \cite{fultonharris}. For each $k\in \N$, we define the $j$-th elementary symmetric 
polynomial $e^k_j(z_1,\ldots,z_k)$, $j=1,\ldots,k$,  in $k$ variables, $z_1,\ldots,z_k$, by:
$$e^k_j=\sum_{i_1<\cdots<i_j}{z_{i_1}z_{i_2}\cdots z_{i_j}}.$$
Equivalently, they can be defined by the following generating function:
$$\sum_{j=0}^{\infty}e^k_j t^j=\prod_{i=1}^k{(1+z_it)}.$$
Obviously, $e^k_j=0$, for $j<0$ or $j>k$, while $e^k_0=1$.

Nice additive basis for homogeneous symmetric polynomials is given by the Schur polynomials. If 
$\lambda=(\lambda_1,\ldots,\lambda_k)$ is a partition such that $\lambda_1\ge\ldots\ge \lambda_k\ge 0$, then the Schur 
polynomial $\pi_{\lambda}(z_1,\ldots,z_k)$ is given by the following expression:
\begin{equation}
\pi_{\lambda}(z_1,\ldots ,z_k)=\frac{|z_i^{\lambda_j+k-j}|}{D},
\label{sur}
\end{equation}
where $D=\prod_{i<j}(z_i-z_j)$, and by $|z_i^{\lambda_j+k-j}|$ we have denoted the determinant of the $k\times k$ 
matrix whose $(i,j)$ entry is equal to $z_i^{\lambda_j+k-j}$. If a $k$-tuple $\lambda=(\lambda_1,\ldots,\lambda_{k})$ doesn't 
satisfy $\lambda_1\ge\cdots\ge\lambda_k\ge 0$, or if $\lambda_{k+1}>0$, then we define $\pi_{\lambda}:=0$.  

From the definition of  $e^k_j$ we have
$$\pi_{1,\ldots,1,0,\ldots,0}=e^k_j,$$
\noindent where there are $j$ units and $k-j$ zeros in the multiindex of $\pi$.

Being the symmetric polynomials, the Schur polynomials can be written as polynomial functions of elementary symmetric 
polynomials. The explicit formula is given by the Giambelli's determinantal formula:
\begin{equation}
\pi_{\lambda}=|e^k_{\mu_i+j-i}|,
\label{djam}
\end{equation}
where $\mu=(\mu_1,\ldots,\mu_l)$, $l=\lambda_1$, is the conjugate partition of $\lambda$, i.e. $\mu_i=\sharp\{j|\lambda_j\ge 
i\}$, and the r.h.s. is the determinant of $l\times l$ matrix.\\

\textbf{Notational convention: From now on, throughout the paper we shall assume that the Schur polynomial is the function of 
elementary symmetric polynomials, and not of the original variables, i.e.:}
$${\pi_{\lambda}(x_1,\ldots,x_k):=|x_{\mu_i+j-i}|.}$$

\vskip 0.3cm

We shall use the following dual form of the Giambelli's determinantal formula: denote by $h_m$ the Schur polynomial
$\pi_{\lambda}$ with $\lambda=(m,0,\ldots,0)$ with $k-1$ zeros (the so-called $m$-th complete symmetric polynomial). Then

\begin{equation}
\pi_{\lambda}=|h_{\lambda_i+j-i}|.
\label{giam}
\end{equation}
For example, we have
$$\pi_{3,3,1}=
\left| \begin{array}{ccc}
h_3 & h_4 & h_5 \\
h_2 & h_3 & h_4 \\
0 & 1 & h_1
\end{array} \right|,$$
where $h_i=\pi_{i,0,0}$.

Moreover, by Giambelli's formula (\ref{djam}), we can express the polynomials $h_m=\pi_{m,0,\ldots,0}$ 
explicitly by:
$$\pi_{m,0,\ldots,0}=f^k_m(x_1,\ldots,x_k).$$
Here the polynomial $f^k_m(x_1,\ldots,x_k)$ is given by the following $m\times m$ determinant:
\begin{equation}
f^k_m(x_1,\ldots ,x_k)= |x_{1+j-i}|,
\label{efovi}
\end{equation}
where by $x_0$ we assume $1$ and $x_i=0$ for $i<0$ or $i>k$. By expanding the determinant (\ref{efovi}) along the first row, 
we obtain:
\begin{equation}
f^k_m=x_1f^k_{m-1}-x_2f^k_{m-2}+\ldots+(-1)^kx_kf^k_{m-k}=\sum_{l=1}^k{(-1)^{l+1}x_lf^k_{m-l}}.
\label{rek}
\end{equation}
Furthermore, by using the Laplace expansion theorem, from (\ref{efovi}) we obtain that for all $p,q\ge 0$ is valid:
\begin{equation}
\pi_{p+q,0,\ldots,0}=\sum_{i=0}^{k-1}{(-1)^i\pi_{p,1,\ldots,1,0,\ldots,0}\pi_{q-i,0,\ldots,0}},
\label{laplas}
\end{equation}
where all multiindices are of length $k$, and there are $i$ $1$'s in the multiindex of the first Schur polynomial on the 
r.h.s. 


\section{Hochschild homology of a bimodule as the homology of a Koszul complex 
of a polynomial ring}\label{3}

Let $R=\C[z_1,\ldots,z_n]$ be the ring of complex polynomials in $n$ variables. We introduce the grading on $R$, by 
defining $\deg z_i=2$, for all $i=1,\ldots,n$. For the partition $i_1,\ldots,i_k$ of $n$ 
(i.e. $n=i_1+\ldots+i_k$) we denote by $R_{i_1,i_2,\ldots,i_k}$ the subring of polynomials invariant under the action of 
the symmetric group $S_{i_1} \times \cdots \times S_{i_k}$. The bimodules obtained as tensor products of such rings over the 
rings of this form are called Soergel bimodules. In this paper we shall focus on the following two $R_{k,l}$-bimodules: 
$B=R_{k,l} \otimes_{R_{k+l}} 
R_{k,l}$ (also denoted as $R_{k,l} \otimes_{{k+l}} R_{k,l}$), and the bimodule $R_{k,l}$ itself. 

In this paper we are interested in describing all graded preserving bimodule maps from $R_{k,l}$ to $R_{k,l} \otimes_{{k+l}} 
R_{k,l}$, as well as in computing the Hochschild homology of the latter bimodule.\\


First of all, as in the case of the Schur polynomials, we shall change our notation for the rings $R_{i_1,\ldots,i_k}$. 
Namely, the polynomial ring $R_{i_1\ldots i_k}$ can be 
represented as the polynomial ring in the variables which are the $i_1$ elementary symmetric polynomials in the first $i_1$ 
variables, the $i_2$ 
elementary symmetric polynomials in the following $i_2$ variables, etc., and the 
$i_k$ elementary symmetric polynomials in the 
last $i_k$ variables. Thus, from now on, we will always work with these ``new" 
variables, i.e. with the elementary symmetric polynomials, because it is more 
convenient for our purposes.

In particular, the ring $R_{k,l}$ is isomorphic to the ring of polynomials in the new variables $x_1,\ldots,x_k$ and 
$y_1,\ldots,y_l$, where each $x_i$ (respectively $y_i$), can be seen as the $i$-th symmetric polynomial in the variables $z_1,
\ldots,z_k$ (respectively $z_{k+1},\ldots,z_{k+l}$). The grading is given by $\deg x_i=2i$ and $\deg y_i=2i$. 
As for the bimodule $R_{k,l} \otimes_{k+l} R_{k,l}$ it is isomorphic to the following quotient module
\begin{equation}
P=\C [x_1,\ldots ,x_k,y_1,\ldots ,y_l,x'_1,\ldots ,x'_k,y'_1,\ldots ,y'_l]/I,
\label{pomkosz}
\end{equation}
where the ideal $I$ is generated by the differences $\Sigma^{x,y}_i-\Sigma^{x',y'}_i$, for all $i=1,\ldots,k+l$, of the 
symmetric polynomials in the variables $x$ and $y$:
\begin{equation}
\Sigma^{x,y}_i=\sum_{j=0}^i {x_j y_{i-j}}.
\label{sigma}
\end{equation}
Here and further on, we assume that $x_0=y_0=1$, $x_i=y_i=0$ for $i<0$, $x_i=0$ for $i>k$ and $y_i=0$ for $i>l$. 
The multiplication by variables $x_i$ and $y_i$ (respectively $x'_i$ and $y'_i$) on $P$ correspond to the left (respectively, 
right) action of $R_{k,l}$ on $B$, while the quotienting by the ideal $I$ matches precisely the 
tensoring over $R_{k+l}$ of the bimodule $B$. Thus, the isomorphism between $P$ and $B$ is given by 
\begin{equation} B \ni p(x,y)\otimes q(x,y) \longleftrightarrow p(x,y)q(x',y') \in P.\label{isom}\end{equation}

Since $I$ is a homogeneous ideal, $P$ inherits the grading from $\C [x_1,\ldots ,x_k,y_1,\ldots ,y_l,x'_1,\ldots ,x'_k,y'_1,\ldots ,y'_l]$. This grading we call the $q$-grading. Also by $\{i\}$ we denote the upward shift by $i$ in the $q$-grading.\\

The Hochschild homology of a bimodule over the polynomial ring can be obtained 
as the homology of the corresponding Koszul complex of certain polynomial rings in many variables. This idea was explained 
and used  in \cite{kovhoh,nas}. Here we shall briefly describe how to 
extend it to our case. 

 The Hochschild homology of $B$ is isomorphic to the homology of the Koszul complex obtained by tensoring 
$$0\xrightarrow{\quad} P\{2i-1\} {\xrightarrow{\ y_i-y'_i\ }} P \xrightarrow{\quad} 0,\quad i=1,\ldots,l,$$
and
$$0\xrightarrow{\quad} P\{2i-1\} {\xrightarrow{\ x_i-x'_i\ }} P \xrightarrow{\quad} 0,\quad i=1,\ldots,k.$$
We put the the right-most terms in the cohomological degree $0$. Since the ideal generated by $I$ and all differences 
$y_i-y'_i$ contains the differences $x_i-x'_i$, the Koszul complex from above is isomorphic to the Koszul complex when the 
differentials in the second group are $0$. Hence, the essential part is in computing the homology of the Koszul complex 
obtained as the tensor product of 
\begin{equation}0\xrightarrow{\quad} P\{2i-1\} {\xrightarrow{\ y_i-y'_i\ }} P \xrightarrow{\quad} 0,\quad i=1,\ldots,l.\label
{koszul}\end{equation}
We denote this complex by $\mathcal{C}$.\\

Now let's focus on the $R_{k,l}$-linear maps $\Delta: R_{k,l} \to B$. If we 
denote $\Delta(1)=\sum_{i}{p_i\otimes q_i}$, for some $p_i,q_i\in R_{k,l}$, then the fact 
that $\Delta$ is a bimodule map is equivalent to $\sum_{i}{wp_i\otimes q_i}=\sum_{i}{p_i\otimes wq_i}$, for all 
polynomials $w\in R_{k,l}$. Hence, $\Delta(1)=\sum_{i}{p_i\otimes q_i}$ defines a bimodule map, if and only if: 
\begin{equation}
\sum_{i}{y_jp_i\otimes q_i}=\sum_{i}{p_i\otimes y_jq_i},\quad \forall j=1,\ldots,l.
\label{bim}
\end{equation}
We pass to the Koszul picture, i.e. we identify $B$ and $P$, and also $R_{k,l}$ and the ring 
$R'=\C[x_1,\ldots,x_k,y_1,\ldots,y_l]$. 

Then $\Delta(1)\in B$, becomes the polynomial $f=\sum_i{p_iq'_i}\in P$, where $q'_i$ is the 
same as the polynomial $q_i$ just with the variables with primes. Then the condition (\ref{bim}) becomes the following set of 
equalities in $P$:
\begin{equation}
(y_j-y'_j)f =0, \quad \forall j=1,\ldots, l,
\label{bim1}
\end{equation}
i.e.
\begin{equation}
(y_j-y'_j)f\in I,\quad \forall j=1,\ldots, l.
\label{bim2}
\end{equation}

Thus, from (\ref{bim1}) we have obtained that $\Delta(1)$ defines a bimodule map, if and only if, the corresponding $f$
defines a class in the homology group $H_{-l}(\mathcal{C})$ of the complex $\mathcal{C}$ from (\ref{koszul}). Hence, we are 
again interested in computation of this homology.\\ 

The Hochschild homology of $B$, as explained above, is isomorphic to the homology of a complex which is the tensor product of 
the complexes of the form 
$$0\ \xrightarrow{\quad\ } R/I \ {\xrightarrow{\,\,\,p\,\,\,}}\ R/I \xrightarrow{\quad\ } 0,$$
where $R$ is a polynomial ring, $I$ is an ideal and $p\in R$ is a polynomial. Let us explain how 
to compute the homology of one such complex. 

The main part is the computation of the kernel and the cokernel of the mapping above. If $p\in I$, the differential is $0$, 
so we assume that $p\notin I$. The cokernel is easily
computed, and is equal to the quotient ring $R/\langle I,p \rangle$. Now, lets pass to the kernel.
For any polynomial $q\in R$ to be in the kernel we must have $pq\in I$.
The ideals are always finitely generated, say $I$ is generated by $i_1,\ldots,i_n$. For any polynomial $q\in R$ to be a 
cocycle we must have $pq\in I$. Therefore we have to find all solutions to the equation $a_1i_1+\cdots+a_ni_n=pq$.
The solutions are generated by $q_1,\ldots,q_k$, which generate a quotient ideal $Q$.
Then the kernel is given by $Q/I$. However, this is isomorphic to $Q  R/\langle I,p \rangle$, since $pq_s\in I$ 
for all $s=1,\ldots,k$. Thus, both kernel and cokernel are of the form $J R/\langle I,p \rangle$, for some ideal $J$. 

\section{The bimodule map}\label{4}

Denote by $\pi_{i_1,\ldots,i_k}$ the Schur polynomial in variables $x_i$'s, i.e. $\pi_{i_1,\ldots,i_k}(x_1,\ldots,x_k)$,  
and denote by $\pi'_{j_1,\ldots,j_l}$ the Schur 
polynomial in variables $y_i$'s, i.e. $\pi_{j_1,\ldots,j_l}(y_1,\ldots,y_l)$ (keep in mind our notational conventions 
from Section \ref{Sur-sim}).
Moreover, if $\alpha=(\alpha_1,\ldots,\alpha_k)$ is the partition with $l\ge\alpha_1\ge\cdots\alpha_k\ge 0$, by $\alpha^*$ we 
denote its complementary partition $\alpha^*=(l-\alpha_k,\ldots,l-\alpha_1)$. Finally by $\bar{\alpha}$ we denote the dual 
(conjugate) partition of $\alpha$, i.e. $\bar{\alpha}_j=\sharp\{i|\alpha_i\ge j\}$. Then, we have the following:

\begin{theorem}
The $R_{k,l}$-linear map $\Delta: R_{k,l} \to R_{k,l} \otimes_{k+l} R_{k,l}$ given by:
\begin{equation}\label{glavna}
\Delta(1)=\sum_{\alpha=(\alpha_1,\ldots,\alpha_k), \\ l\ge\alpha_1\ge\cdots\alpha_k\ge 0} {(-1)^{\sum_{i=1}^k{\alpha_i}} 
\pi_{\alpha}\otimes \pi'_{\bar{\alpha^*}}},
\end{equation}
is a bimodule map. Moreover, it is the only (up to nonzero scalar multiplication) bimodule map of degree $2kl$, and there are 
no nonzero bimodule maps of the degree strictly less than $2kl$.
\end{theorem}

\textbf{Proof:}

Without loss of generality, we can assume that $k\ge l$. 

First we shall prove that $\Delta$ {\em is} a bimodule map. As we explained above, it is enough to show that the polynomial
$$f=\sum_{\alpha=(\alpha_1,\ldots,\alpha_k), \\ l\ge\alpha_1\ge\cdots\alpha_k\ge 0} {(-1)^{\sum_{i=1}^k{\alpha_i}} 
\pi_{\alpha}\pi'_{\bar{\alpha^*}}},$$
where now $\pi'$ is the polynomial in variables $y'$'s, satisfies
$$(y_j-y'_j)f\in I,\quad \forall j=1,\ldots,l.$$
Every element  $r\in I$ can be written as $r=\sum_{i=1}^{k+l} {a_i(\Sigma^i_{x,y}-\Sigma^i_{x',y'})}$, for some $a_i\in R$.
Hence we have to show that for $f$ and every $i_0=1,\ldots,l$, there exists $k+l$ polynomials $a^{i_0}_i\in R$, 
$i=1,\ldots,k+l$,
such that 
$$(y_{i_0}-y'_{i_0})f=\sum_{i=1}^{k+l} {a_i^{i_0}(\Sigma^i_{x,y}-\Sigma^i_{x',y'})}.$$
The right hand side can be written as
\begin{eqnarray*}
\sum_{i=1}^{k+l} {a_i^{i_0}(\Sigma^i_{x,y}-\Sigma^i_{x',y'})}&=&\sum_{i=1}^{k+l} {a_i^{i_0}\left(\sum_{j=1}^i 
x_{i-j}(y_j-y'_j)+
\sum_{j=1}^i y'_{i-j}(x_j-x'_j)\right)}=\\ 
&=&\sum_{j=1}^l (y_j-y'_j)\left(\sum_{i=0}^k a_{i+j}^{i_0}x_i \right) + \sum_{j=1}^k (x_j-x'_j)\left(\sum_{i=0}^l 
a_{i+j}^{i_0}y'_i\right).
\end{eqnarray*}
So, we would have that $(y_{i_0}-y'_{i_0})f$ is equal to the last expression if the multiples of all $x_i-x'_i$ and of all 
$y_i-y'_i$ (except $i_0$-th) are equal to zero, while $f$ is equal to the term that multiplies $y_{i_0}-y'_{i_0}$, i.e.
\begin{eqnarray}
&\sum_{i=0}^k a^{i_0}_{i+j}y'_i = 0, \quad j=1,\ldots,k,&\label{prva}\\
&\sum_{i=0}^l a^{i_0}_{i+j}x_i = 0, \quad j=1,\ldots,l,\quad j\ne i_0,& \label{druga}\\
&f=\sum_{i=0}^l a^{i_0}_{i+i_0}x_i.&\label{treca}
\end{eqnarray}
From the first $l$ equations from (\ref{prva}) we can express $a_l^{i_0},a_{l-1}^{i_0},\ldots,a_1^{i_0}$ in terms of the remaining $a_i^{i_0}$'s as:
\begin{equation}\label{pom1}
a_i^{i_0}=(-1)^{l+1-i}\sum_{j=1}^k \pi'_{l+1-i,{\scriptsize{\underbrace{1,1,\ldots,1}_{j-1}}},0,\ldots,0} a_{l+j}^{i_0},\quad i=1,\ldots,l.
\end{equation}
As we said, by $\pi'_{\alpha}$, we have denoted the Schur polynomial in variables $y'_i$ (and consequently
it has the multiindex of length $l$).
Analogously, from (\ref{druga}) and (\ref{treca}) we obtain
\begin{equation}\label{pom2}
a_i^{i_0}=(-1)^{l+1-i}\sum_{j=1}^k \pi_{l+1-i,{\scriptsize{\underbrace{1,1,\ldots,1}_{j-1}}},0,\ldots,0} a_{l+j}^{i_0},\quad i=1,
\ldots,l, i\ne i_0,
\end{equation}
and
\begin{equation}\label{pom3}
a_{i_0}^{i_0}-f=(-1)^{l+1-i_0}\sum_{j=1}^k \pi_{l+1-i_0,{\scriptsize{\underbrace{1,1,\ldots,1}_{j-1}}},0,\ldots,0}
a_{l+j}^{i_0}.
\end{equation}
Here $\pi_{\alpha}$ is the Schur polynomial in variables $x_i$ (and consequently it has the multiindex of length $k$). 
From (\ref{pom1}), (\ref{pom2}) and (\ref{pom3}) we have that 
\begin{equation}
f=(-1)^{l+1-i_0}\sum_{j=1}^k (\pi'_{l+1-i_0,{\scriptsize{\underbrace{1,1,\ldots,1}_{j-1}}},0,\ldots,0}-\pi_{l+1-i_0,
{\scriptsize{\underbrace{1,1,\ldots,1}_{j-1}}},0,\ldots,0}) a_{l+i_0}^{i_0},
\end{equation}
for certain polynomials $a_{l+1}^{i_0},\ldots,a_{k+l}^{i_0}$ that satisfy the following $l-1$ equations:
$$\sum_{j=1}^k 
(\pi'_{l+1-i,{\scriptsize{\underbrace{1,1,\ldots,1}_{j-1}}},0,\ldots,0}-\pi_{l+1-i,{\scriptsize{\underbrace{1,1,\ldots,1}_{
j-1}}},0,\ldots,0}) a_{l+i}^{i_0}=0, \quad i=1,\ldots,l,\quad i\ne i_0.
$$
Moreover, from the last $k-l$ equation of (\ref{prva}), we have $k-l$ more relations that $a_{l+1}^{i_0},\ldots,a_{k+l}^{i_0}$
should satisfy. Hence, we have that these $k$ polynomials should satisfy certain $k-1$ relations:
$$\sum_{j=1}^k c_i^j a_{l+j}^{i_0} = 0, \quad i=1,\ldots,k-1,$$
while 
$$f=\sum_{j=1}^k c_k^j a_{l+j}^{i_0}.$$
A solution to these equations is given by the following determinant:
$$f=\left| 
\begin{array}{ccc}
c_1^1 & \cdots & c_1^k \\
\vdots & \ddots & \vdots \\
c_k^1 & \cdots  & c_k^k 
\end{array} \right|.$$
Finally, from the form of the coefficients $c_i^j$ in our case, we see that the determinant from above satisfies the 
condition for all $i_0=1,\ldots,l$. 
Hence, if we define the following $k \times k$ matrix
\arraycolsep 6pt
\begin{equation}\label{det}M={{
\left[
\begin{array}{c}
\begin{array}{cccc}
\pi'_{l,0,\cdots,0}-\pi_{l,0,\cdots,0}&\pi'_{l,1, \cdots, 0}-\pi_{l ,1, \cdots, 0}&\cdots&\pi'_{l, 1, \cdots ,1}-\pi_{l, 1, 
\cdots, 1}\\
\pi'_{l-1, 0, \cdots ,0}-\pi_{l-1, 0, \cdots ,0}&\pi'_{l-1, 1, \cdots ,0}-\pi_{l-1, 1, \cdots, 0}&\cdots&\pi'_{l-1, 1, \cdots,
 1}-\pi_{l-1, 1, \cdots, 1}\\
\cdots&\cdots&\cdots&\cdots\\
\pi'_{1, 0, \cdots, 0}-\pi_{1, 0, \cdots ,0}&\pi'_{1, 1, \cdots ,0}-\pi_{1 ,1, \cdots ,0}&\cdots&\pi'_{1, 1, \cdots, 1}-\pi_
{1, 1, \cdots, 1}
\end{array}\\
\hline
N
\end{array}\right],}}\end{equation}
where $N$ is the following $k-l$ by $k$ matrix:
$$N=\left[\begin{array}{cccccccc}
1& y'_1&y'_2&\cdots&y'_l&0&0&0\\
0& 1&y'_1&\cdots&y'_{l-1}&y'_l&0&0\\
\ddots&\ddots&\ddots&\ddots&\ddots&\ddots&\ddots&\ddots\\
0&0&0& 1&y'_1&\cdots&y'_{l-1}&y'_l
\end{array}
\right],$$
then we have that $(y_i-y'_i)\det M\in I$ for all $i=1,\ldots,l$. 
So we are left with proving that $\det{M}$ is (up to a sign) equal to $f$.\\

We shall prove the last fact in the case $k=l$, the remaining cases are done completely analogously. In this case, we have 
that the entries of $M(=[c_i^j])$ are given by 
$$c_i^j=\pi'_{l+1-i,{\scriptsize{\underbrace{1,1,\ldots,1}_{j-1}}},0,\ldots,0}-\pi_{l+1-i,{\scriptsize{\underbrace{1,1,
\ldots,1}_{j-1}}},0,\ldots,0},\quad 1\le i,j\le k.$$
and so the whole determinant can be written as the determinant of the difference of two matrices: one with entries only 
$\pi'$'s and the other only with $\pi$'s. 
We use the formula for the determinant of the sum of two matrices. Let $A$ and $B$ be square $k\times k$ matrices (
$A=[\pi'_{l+1-i,{\scriptsize{\underbrace{1,1,\ldots,1}_{j-1}}},0,\ldots,0}]$ and 
$B=[-\pi_{l+1-i,{\scriptsize{\underbrace{1,1,\ldots,1}_{j-1}}},0,\ldots,0}]$). Then the determinant of $A+B$ is given by 
the signed sum of the products of the determinant of every $i\times i$ submatrix 
of $A$ and the determinant of the $(k-i)\times(k-i)$ complementary submatrix of $B$.  More precisely, let $a=(a_1,\ldots,a_i)$ 
and $b=(b_1,\ldots,b_i)$ be two increasing subsequences of the set $\{1,\ldots,k\}$, and let $a'=(a'_1,\ldots,a'_{k-i})$ and 
$b'=(b'_1,\ldots,b'_{k-i})$ be the complementary increasing subsequences, i.e. such that 
$\{a_1,\ldots,a_i,a'_1,\ldots,a'_{k-i}\}$ $=$ $\{b_1,\ldots,b_i,b'_1,\ldots,b'_{k-i}\}$ $=$ $\{1,\ldots,k\}$. Also, we denote 
by $A_a^b$ the submatrix of $A$ formed by the rows $a_1,\ldots,a_i$ and columns $b_1,\ldots,b_i$. Then we have:
\begin{equation}\label{detz}
\det (A+B) = \sum_{i=0}^k\sum_{a,b} (-1)^{\sum_{j=1}^i{a_j}-\sum_{j=1}^i{b_j}} \det (A_a^b) \det (B_{a'}^{b'}),
\end{equation}
where the second sum is over all increasing subsequences $a$ and $b$ of $\{1,\ldots,k\}$ of length $i$.
The sign is equal to the sign of the summand of the product of entries at the positions $(a_1,b_1), (a_2,b_2),\ldots, 
(a_i,b_i)$, $(a'_1,b'_1), (a'_2,b'_2), \ldots, (a'_{k-i},b'_{k-i})$ in the formula for determinant. This, in turn, is equal 
to the product of the signs of the permutations $(a_1,a_2,\ldots,a_i,a'_1,a'_2,\ldots,a'_{k-i})$ and 
$(b_1,b_2,\ldots,b_i,b'_1,b'_2,\ldots,b'_{k-i})$, which are equal to $\sum_{j=1}^i(a_j-j)$, and $\sum_{j=1}^i(b_j-j)$, 
respectively. This gives the sign from (\ref{detz}).

 Since for each $i=0,\ldots,k$ there are 
$k\choose  i$ subsets of the cardinality $i$, there are 
$$\sum_{i=0}^k {k \choose i}^2 = {2k \choose k}$$  
summands in the expression (\ref{detz}). It can be easily seen that this is equal to the number of partitions 
$\alpha=(\alpha_1,\ldots,\alpha_k)$ with $k\ge\alpha_1\ge\cdots\ge\alpha_k\ge 0$. Indeed, by passing to the sequence of 
differences $(k-\alpha_1,\alpha_1-\alpha_2,\cdots,\alpha_{k-1}-\alpha_k,\alpha_k)$, we have that the number of such 
$\alpha$'s is the same as the number of ways to write $k$ as the sum of $k+1$ nonnegative integers, which is 
$${k+(k+1)-1\choose k+1-1}={2k\choose k}.$$
Moreover, by using the Giambelli 
determinant expression (\ref{giam}) for the Schur polynomials, together with (\ref{laplas}), we have that
\begin{equation}\label{br1}
\det (A_a^b) = (-1)^{i(i-1)/2}\pi'_{(k+1-a_1,k+2-a_2,\ldots ,k+i-a_i,i+1-b'_1,i+2-b'_2,\ldots ,k-b'_{k-i})},
\end{equation}
and 
\begin{equation}\label{br2}
\det (B_{a'}^{b'}) = (-1)^{k-i}(-1)^{(k-i)(k-i-1)/2}\pi_{(k+1-a'_1,k+2-a'_2,\ldots ,k+k-i-a'_{k-i},k-i+1-b_1,k-i+2-b_2,\ldots 
,k-b_i)}.
\end{equation}
Moreover, if we denote by $\alpha$ the multiindex of the Schur polynomial in (\ref{br1}), then the multiindex of the Schur 
polynomial from (\ref{br2}) is given by $\bar{\alpha^*}$. Indeed, this follows from the fact that 
$$a'_l=\sharp \{ j | a_j \le j + l-1\} +l=i+l-\sharp\{j| a_j-j\ge l\}, \quad l=1,\ldots,k-i.$$
Conversely, to every partition  $\alpha=(\alpha_1,\ldots,\alpha_k)$ with $k\ge\alpha_1\ge\cdots\ge\alpha_k\ge 0$, correspond 
sequences $a$ and $b$ defined by $a_j=k+j-\alpha_j$ for $j=1,\ldots,i$, and $b_j=i+j-\alpha_{i+j}$, $j=1,\ldots,k-i$, where 
$i=\max\{j|\alpha_j\ge j\}$.

Hence, we have proved that the determinant from (\ref{det}) is equal to: 
\begin{equation}
\label{kkplus1}
\det M=(-1)^{k(k+1)/2}\sum_{\alpha} (-1)^{\sum_{i=1}^k{\alpha_i}}\pi_{\alpha}\pi'_{\bar{\alpha}^{\ast}},
\end{equation}
and so  $\Delta$ given by (\ref{glavna}) defines a bimodule map, as wanted. 

\begin{remark}
The overall multiplication by $(-1)^{k(k+1)/2}$, although unimportant in our result, precisely matches the ``natural" 
definition of the dual Schur polynomial we obtained in \cite{sln} - Section 3, corresponding to the volume form in the 
cohomology ring of the Grassmanian $Gr_{k,2k}$. More precisely, we have that 
$(-1)^{k(k+1)/2}\pi'_{\bar{\alpha}^{\ast}}=\bar{\pi}_{\alpha^{\ast}}'$, where 
$\bar{\pi}_{\beta}'$ is the dual Schur polynomial of ${\pi}_{\beta}'$.  
\end{remark}

Before going to the uniqueness part, we shall introduce some notation needed for the rest of the paper.  Let 
$R=\C[x_1,\ldots,x_k,y_1,\ldots,y_l,x'_1,\ldots,x'_k,y'_1,\ldots,y'_l]$ (and 
so $P=R/I$), $P_1=P/\langle y_1-y'_1 \rangle$ and $R'=R/\langle I,y_1-y'_1,\ldots,y_l-y'_l\rangle=R/\langle 
x_1-x'_1,\ldots,x_k-x'_k\rangle\simeq \C[x_1,\ldots,x_k,y_1,\ldots,y_l]$. Finally, define the  ideals $I_1,\ldots, I_l$ of 
the ring $R'$
 as follows: the ideal $I_i$ is 
generated by  all $(k-l+i)\times (k-l+i)$ minors of the submatrix of $M$ formed by the last $k-l+i$ rows, 
$i=1,\ldots,l$. In particular $I_l=\langle f\rangle$, as we have just proved.

Now, let $$M'=\left[\begin{array}{cc}S&0\\0&I_{k-l}\end{array}\right]M,$$ where 
$$S={\small{\left[\begin{array}{cccccc}1&-x_1&x_2&\cdots&\cdots&(-1)^{l-1}x_{l-1}\\
0&1&-x_1&\cdots&\cdots&(-1)^{l-2}x_{l-2}\\
0&0&1&\cdots&\cdots&\vdots\\
\cdots&\cdots&\cdots&\cdots&\cdots&\cdots\\
0&0&0&\cdots&1&-x_1\\
0&0&0&\cdots&0&1\end{array}\right].}}$$
Denote the entries of the matrix $M'$ by $c'^j_i$. Moreover, the ideals $I_1',\ldots,I_l'$ defined for the matrix $M'$ in 
the same way as $I_1,\ldots,I_l$ for the matrix $M$, obviously satisfy that $I_j=I_j'$, for all $j=1,\ldots,l$.

Now, it is straightforward to see that $P$ is isomorphic to the ring $\C[x_1,\ldots,x_k,y_1,\ldots,y_l,$ $y_1',\ldots,y_l']/I'$, 
where $I'$ is the ideal generated by the polynomials $(y_1-y'_1){c'^i_1}+(y_2-y'_2){c'^i_2}+\ldots+(y_l-y'_l){c'^i_l}$, for 
all $i=1,\ldots,l$. So we have that 
\begin{equation}
\label{izraz1}
P\simeq R'[Y_1,Y_2,\ldots,Y_l]/J,
\end{equation}
where $J$ is the ideal generated by the polynomials $Y_1{c'^i_1}+Y_2{c'^i_2}+\ldots+Y_l{c'^i_l}$, for 
all $i=1,\ldots,l$. New variables $Y_i$ have degrees $2i$, $i=1,\ldots,l$.\\

\textit{Uniqueness:}\quad Observe the following complex:
\begin{equation}
0\rightarrow P\{2\} {\xrightarrow{y_1-y'_1}} P \rightarrow 0.
\label{pomoc1}
\end{equation}
 The homology at the right-most term (denoted $H'_0$) is 
equal to $P_1=P/\langle y_1-y'_1 \rangle$, while the homology at the the left-most term (denoted $H'_1$), contains the ideal 
$f P$, as we have proved above. Moreover, since $P$ is of the form $R/I$, and we have proved that $(y_i-y'_i)f \in I$, for 
all $i=1,\ldots,l$, we have that $f P$ is isomorphic to $f P/\langle y_i-y'_i \rangle$, where the quotient is taken over all 
differences $y_i-y'_i$, for all $i=1,\ldots,l$. However, $R/\langle I,y_i-y'_i \rangle$ is isomorphic to the ring 
$R'$, 
and so $f P$ is isomorphic to $f R'$.

Since the differential in (\ref{pomoc1}) is grading preserving, the graded Euler 
characteristic of the complex (\ref{pomoc1}) satisfies:
\begin{equation} (1-q^{2}) \qdim P =  \qdim H'_0 - q^{2}\qdim H'_1,\label{11}\end{equation} 
and so 
\begin{equation}
\label{nova}
\qdim H'_1=q^{-2}(\qdim P_1-(1-q^{2})\qdim P).
\end{equation}
Below we shall prove that 
\begin{equation}
\label{nova1}
(1-q^{2})\qdim P=\qdim P_1-q^2 q^{2kl}\qdim R'.
\end{equation}

Then from (\ref{nova}) and (\ref{nova1}), we have 
$$\qdim H'_1=q^{2kl}\qdim R'.$$
Since $\deg f =2kl$ we have that $\qdim (fR') = \qdim H_1'$, and since we have 
proved that $fR' \subset H_1'$, we have that $H_1'=fR'$. Hence all $g$ that satisfy $(y_1-y'_1)g \in I$, are  multiples of 
$f$, which proves the uniqueness part.\\

Thus we are left with proving (\ref{nova1}). We shall  use the expression (\ref{izraz1}) for the ring 
$P$. By eliminating variables $Y_l,Y_{l-1},\ldots,Y_1$ we have that 
\begin{eqnarray*}
P&\simeq & R'[Y_1,\ldots,Y_{l-1}]/(J\cap R'[Y_1,\ldots,Y_{l-1}])\oplus Y_lR'[Y_1,\ldots,Y_{l}]/(Y_l I_1)\\
&\simeq &\cdots \quad\simeq\quad R'\oplus \oplus_{i=1}^lY_i R'[Y_1,\ldots,Y_{i}]/(Y_i I_{l+1-i}).\end{eqnarray*}
Thus, we obtained that 
\begin{equation}\qdim P=\qdim R'+\sum_{i=1}^l{\frac{q^{2i}}{\prod_{j=1}^i{(1-q^{2j})}}}\qdim R'/I_{l+1-i}.\label{nova2}\end{equation}
By using the expression (\ref{izraz1}),  we have that $P_1$ is isomorphic to the ring $R'[Y_2,\ldots,Y_l]/J_1$, where $J_1$ is the ideal 
generated by the following $l$ polynomials: $Y_2c'^i_2+Y_3c'^i_3+\ldots+Y_lc'^i_l$, $i=1,\ldots,l$. Then, as above, we have 
\begin{eqnarray*}
P_1&\simeq & R'[Y_2,\ldots,Y_{l-1}]/(J_1\cap R'[Y_2,\ldots,Y_{l-1}])\oplus Y_lR'[Y_2,\ldots,Y_{l}]/(Y_l I_1)\\
&\simeq &\cdots \quad\simeq\quad R'\oplus \oplus_{i=2}^lY_i R'[Y_2,\ldots,Y_{i}]/(Y_i I_{l+1-i}),\end{eqnarray*}
and so 
 \begin{equation}
\qdim P_1= \qdim R' +\sum_{i=2}^l {{\frac{q^{2i}}{\prod_{j=2}^i {(1-q^{2j})}}}\qdim R'/I_{l+1-i}}.
\label{pe1}
\end{equation}

Hence from (\ref{nova2}) and (\ref{pe1}), we have 
\begin{equation}
\label{nova3}
\qdim P=\qdim R'+\frac{q^2}{1-q^2}\qdim R'/I_l+\frac{1}{1-q^2}(\qdim P_1-\qdim R'),
\end{equation}
and since $I_l=\langle f\rangle$, we obtain
\begin{equation}
\label{nova5}
(1-q^2)\qdim P=\qdim P_1-{q^2}\qdim I_l=\qdim P_1-{q^2}q^{2kl}\qdim R',
\end{equation}
which gives the wanted formula (\ref{nova1}). 
\kraj

\begin{remark}
As in the previous proof, we can obtain more general result. Namely, let $P_t=P/\langle y_1-y'_1,\ldots,y_t-y'_t \rangle$, 
$t=0,\ldots,l$. Obviously, $P_0=P$ and $P_l=R'$. Then, by using the same expression (\ref{izraz1}), we obtain that for every 
$t=0,\ldots,l-1$:
\begin{eqnarray*}
\qdim P_t &=& \qdim R' + \frac{q^{2t}}{1-q^{2t}}\qdim R'/I_{l-t} + \frac{1}{1-q^{2t}}(\qdim P_{t+1} -\qdim R')=\\
&=& \frac{1}{1-q^{2t}} \qdim P_{t+1} - \frac{q^{2t}}{1-q^{2t}} \qdim I_{l-t},
\end{eqnarray*}
and so:
\begin{equation}
\label{vazna}
q^{2t}\qdim I_{l-t} =\qdim P_{t+1}-(1-q^{2t})\qdim P_t,\quad t=0,\ldots, l-1.
\end{equation}
\end{remark}

\section{Hochschild homology of the bimodule $B$}\label{5}

Let $B=R_{k,l}\otimes_{k+l}R_{k,l}$, and let $P$ be the corresponding presentation (\ref{isom}).  The aim is to compute the 
homology groups $H_{-i}$, for $i=0,\ldots,l$ of the Koszul complex obtained by tensoring
$$0\rightarrow P\{2i-1\} {\xrightarrow{y_i-y'_i}} P \rightarrow 0,$$
for $i=1,\ldots,l$.
\begin{theorem}
The homology groups $H_{-i}$, for $i=0,\ldots,l$, are given by 
$$H_{-i}\simeq \bigoplus_{1\le \alpha_1 <\cdots < \alpha_i\le l} I_{l+1-\alpha_1} \{2\sum_{m=1}^i a_m - i\},\quad l\ge i>0,$$
$$H_0 \simeq R'.$$
\end{theorem}

\textbf{Proof:}

To compute the homology of the Koszul complex, we first take the homology with respect to the first differential  
(i.e. strip-off the contractible summand), then take the homology induced by  the second differential, and so on.
Denote the homology obtained in this way after taking the homology with respect to the $j$-th differential by $H^j$.
We have that $H^j_i$ is trivial for $i<-j$ or $i>0$. We shall prove by induction on $j$ that 
\begin{equation}
H^j_{-i}\simeq \bigoplus_{1\le \alpha_1 <\cdots < \alpha_i\le j} I_{l+1-\alpha_1} \{2\sum_{m=1}^i a_m - i\}, \quad j\ge i>0,
\label{homj}
\end{equation}
and
\begin{equation}
H^j_0 \simeq P_j,
\label{hom0}
\end{equation}
where $P_j$ is the quotient of the polynomial ring $P$ by the ideal generated by $y_i-y'_i$ for $1\le i\le j$.

For $j=1$, the homology $H^1$ is the homology of the complex
$$0\rightarrow P\{1\} {\xrightarrow{y_1-y'_1}} P \rightarrow 0.$$
For this homology we have obtained in the previous section that $H^1_{-1}=I_l\{1\}$, while $H^1_0=P_1$, in accordance with 
(\ref{homj}) and (\ref{hom0}).

Now, suppose that the homology $H^j$ for some $1\le j<l$ is given by (\ref{homj}) and (\ref{hom0}), and lets compute the 
homology $H^{j+1}$. It is given by the homology of the complex
$$0\rightarrow H^j\{2j+1\} {\xrightarrow{y_{j+1}-y'_{j+1}}} H^j \rightarrow 0.$$
The differential is zero on all $I_i$ (since they are the ideals of the ring $R'$), and so the differential is nontrivial 
only on 
$$0\rightarrow P_j\{2j+1\} {\xrightarrow{y_{j+1}-y'_{j+1}}} P_j \rightarrow 0.$$
This homology can be computed completely analogously as in Theorem 1 (see also Remark 2 - formula (\ref{vazna})), and the 
homology at the left-most term is isomorphic to $I_{l-j}\{2j+1\}$, while the homology at the right-most term is $P_{j+1}$. 
Hence, altogether we have that
\begin{eqnarray*}
H_{-i}^{j+1}&\simeq &H_{-i}^j \oplus H_{-i+1}^j\{2j+1\}, \quad j+1\ge i\ge 2,\\
H_{-1}^{j+1}&\simeq &H_{-1}^j \oplus I_{l-j}\{2j+1\},\\
H_0^{j+1}&\simeq & P_{j+1},
\end{eqnarray*}
which together with (\ref{homj}) gives the wanted formula. \kraj

\begin{remark}
As explained at the end of the section \ref{3}, all homology groups $H_{-i}$ are isomorphic to the direct sums of certain 
ideal over the ring $P/\langle y_1-y'_1,\ldots,y_l-y'_l\rangle \simeq R'$. Also, they can be written by the following formal 
expression  
$$\sum_{i=0}^l{(-1)^it^{-i}H_{-i}}=(1-t^{-1}q^{2l-1}I_1)(1-t^{-1}q^{2l-3}I_2)\cdots(1-t^{-1}qI_l),$$
where by the product $I_{\alpha_1}I_{\alpha_2}\cdots I_{\alpha_j}$ we mean $I_{\max\{\alpha_1,\ldots,\alpha_j\}}$. 
\end{remark}

As a corollary of the previous theorem, we obtain the categorification of one digon move axiom of the calculus of the colored 
HOMFLY-PT polynomial (axiom A2 from \cite{nas}). 
\begin{corollary}
$$\sum_{i=0}^l{(-1)^i t^i\qdim H_{-i}}=\qdim R' \prod_{i=1}^l (1-t^{-1}q^{2k+2i-1}). $$
\end{corollary}
\textbf{Proof:}

First of all, the right hand side of the equation above, can be expanded as
\begin{equation}\label{rhs}
\prod_{i=1}^l (1-t^{-1}q^{2k+2i-1})=\sum_{i=0}^l(-1)^i t^{-i} q^{i(i+2k)}\left[l\atop i\right].
\end{equation} 
Here  we use the following notation for quantum integers: 
\begin{eqnarray*}
[n]&=&\frac{1-q^{2n}}{1-q^2}=1+q^2+\ldots +q^{2(n-1)},\\
{[n]!}&=& [n][n-1]\cdots[2][1],\\
\left[n \atop m\right] &=& \frac{[n]!}{[m]![n-m]!}.
\end{eqnarray*}
On the other hand, by Theorem 2 (Remark 3),  we have 
\begin{eqnarray}
&&\sum_{i=0}^l{(-1)^i t^i\qdim H_{-i}}=\qdim R'-\nonumber\\
&&-\sum_{j=1}^l t^{-1}q^{2l+1-2j}(1-t^{-1}q^{2l+3-2j})(1-t^{-1}q^{2l+5-2j})\cdots(1-t^{-1}q^{2l-1})\qdim I_j.\label{novo7}\end{eqnarray}
Now, we are left with computing $\qdim{I_j}$, $j=1,\ldots,l$. Since  
$I_j$ is generated by the minors 
of the matrix whose entries constitute a regular sequence, there is a natural free graded resolution of it, from which we obtain 
\begin{equation}\label{nova4}\qdim I_j= \qdim R' \sum_{i=j}^l {(-1)^{i-j}}q^{i(i+1+2(k-l))+j(j-1)}\left[i-1 \atop i-j\right]\left[l \atop i\right], 
\quad j=1,\ldots,l.
\end{equation}
By replacing (\ref{nova4}) in (\ref{novo7}),  and since  $$\prod_{i=1}^{j-1}(1-t^{-1}q^{2l+1-2j+2i})=\sum_{s=0}^{j-1}(-1)^s t^{-s} q^{(2l+1-2j)s} q^{s(s+1)}\left[{j-1}\atop{s}\right],$$ we have 
\begin{eqnarray}&&\sum_{i=0}^l{(-1)^i t^i\qdim H_{-i}}=\qdim R'\cdot\label{n9}\\
&&\cdot(1-\sum_{j=1}^l\sum_{s=0}^{j-1}\sum_{i=j}^l{(-1)^{i-j+s}}t^{-s-1}q^{(2l+1-2j+s)(s+1)}q^{i(i+1+2(k-l))+j(j-1)}\left[{j-1}\atop s\right]\left[{i-1}\atop i-j\right]\left[{l}\atop i\right]).\nonumber\end{eqnarray}

By exchanging the order of summations, the last triple sum becomes
\begin{equation}\sum_{i=1}^l(-1)^i q^{i(i+1+2(k-l))}\left[{l}\atop i\right]\left(\sum_{s=0}^{i-1}(-1)^s q^{(2l+1+s)(s+1)} t^{-s-1}\left(
\sum_{j=s+1}^{i}(-1)^j q^{-2j(s+1)} q^{j(j-1)}\left[{i-1}\atop i-j\right]\left[{j-1}\atop s\right]\right)\right).\label{novo8}\end{equation}
Since $$\left[{i-1}\atop i-j\right]\left[{j-1}\atop s\right]=\frac{[i-1]!}{[i-j]![j-1]!} \frac{[j-1]!}{[s]![j-1-s]!}=\left[{i-1}\atop s\right]\left[{i-1-s}\atop j-1-s\right],$$
the innermost sum in  (\ref{novo8}) becomes
$$\left[{i-1}\atop s\right]\sum_{j=s+1}^{i}(-1)^{j}q^{j(j-2s-3)}\left[{i-1-s}\atop {j-1-s}\right]=\left[{i-1}\atop s\right]\sum_{x=0}^{i-1-s}(-1)^{x+s+1}q^{(x+s+1)(x-s-2)}\left[{i-1-s}\atop x\right]=$$
$$(-1)^{s+1}q^{-s^2-3s-2}\left[{i-1}\atop s\right]\sum_{x=0}^{i-1-s}(-1)^{x}q^{x(x-1)}\left[{i-1-s}\atop x\right]=
(-1)^{s+1}q^{-s^2-3s-2}\left[{i-1}\atop s\right]\delta_{i-1-s,0},$$
where $\delta_{a,b}$ is the Kronecker delta. By replacing this in (\ref{novo8}) we obtain
$$-\sum_{i=1}^l(-1)^i q^{i(i+1+2(k-l))}\left[{l}\atop i\right]\left(\sum_{s=0}^{i-1}q^{(2l-1)(s+1)} t^{-s-1}\left[{i-1}\atop s\right]\delta_{i-1-s,0}\right)=$$
$$=-\sum_{i=1}^l(-1)^i q^{i(i+2k)}t^{-i}\left[{l}\atop i\right].$$
Finally, by replacing the obtained value of the triple sum in (\ref{n9}), we obtain 
$$\sum_{i=0}^l{(-1)^i t^i\qdim H_{-i}}=\qdim R'\sum_{i=0}^l(-1)^i q^{i(i+2k)}t^{-i}\left[{l}\atop i\right],$$
which together with (\ref{rhs}) finishes the proof.

\kraj

\end{document}